
\magnification=1100
\baselineskip=14pt

\def\nin{\noindent}

  \def\FF{{\cal F}}  \def\GG{{\cal G}}

      \def\PPP{{\bf P}}

\def\IR{{{\rm I}\!{\rm R}}}

\def\ref#1{{\rm [}{\bf #1}{\rm ]}}   
\def\nref#1#2{{\rm [}{\bf #1}{\rm ;\ #2]}}

\def\comp{{\leavevmode
     \raise.2ex\hbox{${\scriptstyle\mathchar"020E}$}}}



\outer\def\proclaim#1{\medbreak\noindent\bf\ignorespaces
   #1\unskip.\enspace\sl\ignorespaces}
\outer\def\endproclaim{\par\ifdim\lastskip<\medskipamount\removelastskip
   \penalty 55 \fi\medskip\rm}

\def\rect#1#2#3{\raise .1ex\vbox{\hrule height.#3pt
   \hbox{\vrule width.#3pt height#2pt \kern#1pt\vrule width.#3pt}
        \hrule height.#3pt}}

\def\qed{$\hskip 5pt\rect364$}

\def\ee{{\bf e}} 

\def\DS{DS53}
\def\FIT{F85}
\def\IMK{IM74}
\def\ITO{I70}
\def\KNT{K82}
\def\KNG{K02}
\def\LEH{L77}

\centerline{\bf A Proof of Lehoczky's Theorem on Drawdowns}
\medskip
\centerline{P.J. Fitzsimmons}
\medskip
\centerline{Department of Mathematics}
\centerline{University of California, San Diego}
\centerline{Department of MathematicsLa Jolla, CA 92093--0112}
\centerline{\tt pfitzsim@ucsd.edu}
\bigskip
\bigskip

\nin{\bf 1. Introduction.
}
 Let $X=(X_t)_{t\ge 0}$ be a regular diffusion process on an interval
$E\subset\IR$. Let $M_t:=\max_{0\le u\le t}X_u$ denote the past maximum process
of $X$, and for $\delta>0$ define the ``drawdown time" $\tau=\tau_\delta$ by
$$
\tau=\inf\{t>0: X_t=M_t-\delta\}.
\leqno(1.1)
$$
Our goal here is to use the excursion theory developed in \ref{\FIT} to prove (the general form of) a theorem of J.~Lehoczky  concerning the joint distribution of $\tau$ and $M_\tau$.

Before stating the result we introduce the necessary notation. Let $A=\inf E$, $B=\sup E$, and write $E^\circ=]A,B[$. We assume throughout the
paper that $B\notin E$, and that $A\in E$ if and only if $A$ is a regular
boundary point which is not a trap for $X$. These assumptions
imply that the transition kernels of $X$ are absolutely continuous with respect
to the speed measure $m$ (recalled below). See \S 4.11 of \ref{\IMK}.

The process $X$ is realized as the coordinate process on the space $\Omega$ of  paths $\omega\colon[0,+\infty[\to E\cup\{\Delta\}$ which are absorbed in the
cemetery point $\Delta\notin E$ at time $\zeta(\omega)$, and which are
continuous on $[0,\zeta(\omega)[$.   The $\sigma$-fields
$\FF$ and $\FF_t$ ($t\ge 0$) are the usual Markovian completions of
$\FF^\circ=\sigma\{X_u:u\ge 0\}$ and $\FF^\circ_t=\sigma\{X_u:0\le u\le t\}$
respectively. The law $\PPP^x$ on $(\Omega,\FF^\circ)$ corresponds to $X$ started
at $x\in E$.

We assume that $X$ admits no killing in $E^\circ$, and let $S$ and $m$ denote the scale function and speed measure of $X$. 
The infinitesimal  generator
$\GG$ of $X$ has the form
$$
\GG u(x)\cdot m(dx)=du^+(x),\qquad x\in E^\circ,
\leqno(1.2)
$$
for $u\in D(\GG)$, the domain of $\GG$. Here and elsewhere $u^+$ denotes the (right-hand) 
scale derivative:
$$
u^+(x)=\lim_{y\downarrow x}{u(y)-u(x)\over S(y)-S(x)}.
\leqno(1.3)
$$
Likewise, $u^-$ denotes the left scale derivative.  

As a regular diffusion, $X$ admits local time; this is a jointly continuous (adapted) process $(L^y_t:t\ge 0, y\in E)$ normalized to be occupation density relative to $m$; that is
$$
\int_0^t f(X_s)\,ds =\int_E f(x)L^x_t\,m(dx),\qquad\forall x\in E, t\ge 0, 
\leqno(1.4)
$$
for all bounded continuous $f$, almost surely. 
For a fixed level $y\in E$, the local time can be used to normalize the 
It\^o excursion law \ref{\ITO}, for excursions from level $y$, as follows. Let
$G(y)$ denote the (random) set of left-hand endpoints (in $]0,\zeta[$) of
intervals contiguous to the level set $\{t>0:X_t=y\}$. Define the hitting time
$T_y$ by
$$
T_y=\inf\{t>0:X_t=y\}\qquad (\inf\emptyset=+\infty).
\leqno(1.5)
$$
The It\^o excursion law $n_y$ is determined by the identity
$$
\PPP^x\sum_{u\in G(y)}Z_u\,F(\ee^u) =\PPP^x\left(\int_0^\infty
Z_u\,dL^y_u\right)\cdot n_y(F),
\leqno(1.6)
$$
where $x\in E$, $\ee^u$ is the excursion from $y$ starting at time $u\in G$
$$
{\ee}^u_t=\cases{ X_{u+t},&$0\le t<T_y-u$,\cr \Delta,&$t\ge T_y-t$.\cr}
\leqno(1.?)
$$
$F\in p\FF^\circ$, and $Z\ge 0$ is an $(\FF_t)$-optional
process. Under $n_y$ the coordinate process $(X_t:t> 0)$ is strongly
Markovian with semigroup $(Q^y_t)$ given by
$$
Q^y_t(x,f)=\PPP^x(f\comp X_t;t<T_y).
\leqno(1.7)
$$

Finally, the point process of excursions below the maximum is defined as
follows. For $t\ge 0$ set
$$
\eqalign{
H(\omega)&=\{u>0:X_u(\omega)=M_u(\omega)\};\cr
R_t(\omega)&=\inf\{u>0:u+t\in H(\omega)\};\cr
G(\omega)&=\{u>0: u<\zeta(w), R_{u-}(\omega)=0<R_u(\omega)\}.\cr}
\leqno(1.8)
$$
\medskip
\noindent 
Thus $G$ is the random set of left-hand endpoints of intervals
contiguous to the random set $H$. As before, for each  $u\in G$ we have an excursion ${\ee}^u$
defined by
$$
{\ee}^u_t=\cases{ X_{u+t},&$0\le t<R_u$,\cr \Delta,&$t\ge R_u$.\cr}
\leqno(1.9)
$$
The point process $\Pi=\{(X_u,{\ee}^u):u\in G\}$ admits a L\'evy system as detailed in the following proposition. In effect, the excursions below the running maximum, when indexed by their levels, form a Poisson point process in $E\times\Omega$, with intensity $(dy,d\omega)\mapsto dS(y)\cdot n^\downarrow_y(d\omega)$. 
Define
a continuous increasing adapted process $C=(C_t:t\ge 0)$ by
$$
C_t=\cases{S(M_t)-S(M_0),&if $t<\zeta$,\cr C_{\zeta-}&if $t\ge\zeta$.\cr}
\leqno(1.10)
$$

\proclaim{(1.11) Proposition} For $Z\ge 0$ and $(\FF_t)$-optional, and $F\in
p\FF^\circ$,
$$
\eqalign{
\PPP^x\sum_{u\in G} Z_u\, F({\ee}^u)&=\PPP^x\int_0^\infty
Z_u\,n^\downarrow_{X_u}(F)\,dC_u\cr
&=\PPP^x\int_A^x Z_{T_y} 1_{\{T_y<+\infty\}} n_y^\downarrow(F)\,dS(y),\cr}
\leqno(1.12)
$$
where $n^\downarrow_y$ denotes the restriction of $n_y$ to
$\{\omega:\omega(t) < y,\forall t\in]0,\zeta(\omega)[\}$.
\endproclaim\bigskip

\nin{\bf (1.13) Remark.} The second equality in (1.12) follows from the first
by the change of variable $u=T_y$. 
\medskip

We will also need (see \S 4.6 of \ref{\IMK}) the Laplace transform  of $T_y$: 
$$
\PPP^x(e^{-\alpha T_y})=\cases{
g^\alpha_1(x)/g^\alpha_1(y),&$x\le y$,\cr g_2^\alpha(x)/g^\alpha_2(y),&$x\ge
y$,\cr}
\leqno(1.14)
$$
where for fixed  $\alpha>0$, $g_1^\alpha$ and $g_2^\alpha$ are strictly positive, linearly independent solutions  of
$$
\GG g(x)=\alpha g(x),\qquad x\in E^\circ;
\leqno(1.15)
$$
$g^\alpha_1$ (resp.\ $g^\alpha_2$) is an increasing (resp.\ decreasing)
solution of (1.15) which also satisfies the appropriate boundary condition at
$A$ (resp.\ $B$). Both $g^\alpha_1$ and $g^\alpha_2$ are uniquely determined up
to a positive multiple. As a rule we drop the superscript $\alpha$, writing
simply $g_1$ and $g_2$. 
 
Define
$$
b(z):=n^\downarrow_z\left(e^{-\alpha T_{z-\delta}};T_{z-\delta}<\zeta\right)
\leqno(1.16)
$$
and
$$ 
c(z):=n^\downarrow_z\left(1-e^{-\alpha\zeta}; T_{z-\delta}>\zeta\right).
\leqno(1.17)
$$
(More explicit expression for $b$ and $c$, in terms of the solutions $g_1$ and $g_2$,  will be displayed below in the course of proving the Theorem.)

\proclaim{(1.18) Theorem}
The joint Laplace transform of $M_\tau$ and $\tau$  is given by
$$
\PPP^x[\exp(-\alpha\tau-\beta  M_\tau)] =
\int_x^B \exp\left(-\beta y-\int_x^y b(z)\,dS(z)\right)\,c(y)\,dS(y),\qquad x\in E.
\leqno(1.19)
$$
In particular, 
$$
\PPP^x[M_\tau >y] =\exp\left(-\int_x^y {dS(z)\over S(z)-S(z-\delta)}\right),\qquad x\le y<B.
\leqno(1.20)
$$ 
\endproclaim
\bigskip

\nin{\bf 2. Proofs.}

\proclaim{(2.1) Lemma} For $\alpha>0$ and $\delta>0$,
$$
n^\downarrow_y\left(e^{-T_{y-\delta}}; T_{y-\delta}<\zeta\right)
={g_1^-(y)g_2(y)-g_1(y)g^-_2(y)\over g_1(y)g_2(y-\delta)-g_1(y-\delta)g_2(y)}.
\leqno(2.2)
$$
\endproclaim

\nin{\sl Proof.} By a result apparently due originally to \ref{\DS} (see also \nref{IMK}{\S 4.10}), for $y-\delta<y-\epsilon<y$, we have
$$
\PPP^{y-\epsilon}\left(e^{-\alpha T_{y-\delta}}; T_{y-\delta}<T_y\right)={g_1(y)g_2(y-\epsilon)-g_1(y-\epsilon)g_2(y)\over g_1(y)g_2(y-\delta)-g_1(y-\delta)g_2(y)}.
\leqno(2.3)
$$
Subtract and add $g_1(y)g_2(y)$ in the numerator on the right, and then divide both sides by $S(y)-S(y-\epsilon)$ and send $\epsilon\downarrow 0$. The limit on the left is
(as is well known, and easily deduced) $n^\downarrow_y\left(\exp(-\alpha T_{y-\delta});T_{y-\delta}<\zeta\right)$. The left-hand scale derivatives of $g_1$ and $g_2$ exist because both functions are in the domain of the generator, so the limit of the right side is as indicated.\qed
\medskip

Because $g_1$ and $g_2$ are linearly independent, any other linearly independent pair of solutions of (1.15) may be used in (2.2) instead, and likewise in (2.5) below.

Essentially the same computation yields
\proclaim{(2.4) Lemma}
$$
n^\downarrow_y\left(1-e^{-\alpha\zeta}; T_{y-\delta}>\zeta\right)={g_1^-(y)g_2(y-\delta)-g_1(y-\delta)g^-_2(y)\over g_1(y)g_2(y-\delta)-g_1(y-\delta)g_2(y)}.
\leqno(2.5)
$$
\endproclaim

The formulas presented in these two lemmas are related to the spectral decomposition of the Laplace transform of the hitting times of $X$ that is treated in \ref{\KNT}.
\medskip
\nin{\sl Proof of the Theorem.} We start with (1.20). Observe that when $X_0=x$, we have $\{M_\tau>y\}=\{N_{x,y}=0\}$, where 
$$
N_{x,y}:=\{(z,\ee)\in\Pi: x<z\le y, T_{z-\delta}(\ee)<\zeta(\ee)\}.
\leqno(2.6)
$$
By Proposition (1.11), under $\PPP^x$ the random variable $N_{x,y}$ has the Poisson distribution with mean value
$$
\int_{]x,y]} n^\downarrow_z[T_{z-\delta}<\zeta]\,dS(z).
\leqno(2.7)
$$
Notice that the $\alpha=0$ version of (2.3) is the well known
$$
\PPP^{z-\epsilon}(T_{z-\delta}<T_z)={S(z)-S(z-\epsilon)\over S(z)-S(z-\delta)}
\leqno(2.8)
$$
and so 
$$
n^\downarrow_z(T_{z-\delta}<\zeta) ={1\over S(z)-S(z-\delta)},
\leqno(2.9)
$$
as in the proof of (2.2). Formula  (1.20) now follows from a well-known property of the Poisson distribution.

We now turn to the proof of (1.19). Taking a cue from \ref{\LEH}, we compute
$$
\PPP^x[e^{-\alpha \tau}\,|\, M_\tau=y]
\leqno(2.10)
$$
for $y>x$.
Notice that, given that $M_\tau=y$,  $\tau$ is the sum of the ``run-up time'' to level $y$
$$
\rho:=\sum_{(z,\ee)\in\Pi} \zeta(\ee)\cdot1_{x<z\le y, T_{z-\delta}(\ee)>\zeta(\ee)}
\leqno(2.11)
$$
(because the Lebesgue measure of $H$ is $0$, a.s.) and a random variable $\sigma$ that is independent of $\rho$ and which has the same distribution as $T_{y-\delta}$ under $n^\downarrow_y(\,\cdot\, | T_{y-\delta}<\zeta)$. 
By the Poisson master formula \nref{\KNG}{(3.6)}
$$
\PPP^x[\exp(-\alpha\rho)]=\exp\left(-\int_{]x,y]} n^\downarrow_z(1-e^{-\alpha \zeta})\,dS(z)\right)
\leqno(2.12)
$$
Combining this with Lemmas (2.1) and (2.4), we obtain
$$
\eqalign{
\PPP^x[e^{-\alpha \tau}\,|\, M_\tau=y]
&=\exp\left(-\int_{]x,y]} n^\downarrow_z(1-e^{-\alpha \zeta},T_{z-\delta}>\zeta)\,dS(z)\right)\cdot n^\downarrow_y(e^{-\alpha T_{y-\delta}} |T_{y-\delta}<\zeta)\cr
&=\exp\left(-\int_x^y c(z)\,dS(z)\right)\cdot b(y)\cdot (S(y)-S(y-\delta)).\cr
}
\leqno(2.13)
$$
In view of (1.20), which shows that
$$
\PPP^x[M_\tau\in dy] = {dS(y)\over S(y)-S(y-\delta)},\qquad y>x,
$$
the Theorem is proved.\qed
 
\bigskip
\centerline{\bf References}
\frenchspacing
\medskip

\item{[\DS]}
Darling, D.A., Siegert, A.J.F.: The first passage problem for a continuous Markov process, {\it Ann. Math. Statistics} {\bf 24} (1953) 624--639.
\smallskip
 
\item{[\FIT]}
Fitzsimmons, P.J.: Excursions above the minimum, \hfil\break {\tt https://arxiv.org/abs/1308.5189}, 1985.
\smallskip

\item{[\IMK]}
It\^o, K. and McKean, H.P.: {\it Diffusion Processes and their Sample Paths}, (Second
printing, corrected.) Springer-Verlag, Berlin-New York, 1974.
\smallskip

\item{[\ITO]}
It\^o, K.: Poisson point processes attached to Markov processes, {\it  Proceedings of
the Sixth Berkeley Symposium on Mathematical Statistics and Probability\/} (Univ.
California, Berkeley, Calif., 1970/1971), Vol. III, pp. 225Ð239. UC Press, Berkeley, 1972.
\smallskip

\item{[\KNT]}
Kent, J.T.: The spectral decomposition of a diffusion hitting time, {\it Ann. Probab.} {\bf 10} (1982) 207--219.
\smallskip

\item{[\KNG]}
Kingman, J.F.C.: {\it Poisson Processes},  Oxford Studies in Probability, Clarendon Press, Oxford, 2002. 
\smallskip

\item{[\LEH]} 
Lehoczky, J.: Formula for stopped diffusion processes with stopping times based on the maximum, {\it Ann. Probab.} {\bf  5} (1977) 601--607.
\smallskip

\end